\RequirePackage{fix-cm}
\pdfoutput=1

\documentclass[a4paper]{article}

\usepackage{listings}
\usepackage{graphicx}
\usepackage{amsmath}
\usepackage{float}
\usepackage{blindtext}
\usepackage{color}
\usepackage[table]{xcolor}
\usepackage{subfig}
\usepackage{tabu}
\usepackage{booktabs}
\usepackage{rotating}
\usepackage{array,multirow}
\usepackage{amssymb}
\usepackage{quoting}
\usepackage{algorithm}
\usepackage[noend]{algpseudocode}
\usepackage[affil-it]{authblk}
\quotingsetup{font=small}

\newcommand{\redmark}[1]{{#1}}
\newcommand{\bluemark}[1]{{#1}}

\newcommand{\ignore}[1]{}
\newcommand{\citepar}[1]{\cite{#1}}

\begin{document}

\title{A general and scalable matheuristic for fleet design
}

\author{Francesco Bertoli\thanks{francesco.bertoli@data61.csiro.au}}
\author{\,Philip Kilby} 
\author{Tommaso Urli}
\affil{108 North Road, Acton ACT 2601, Australia}


\maketitle

\begin{abstract}
  We look at the problem of choosing a fleet of vehicles to carry out
  delivery tasks across a very long time horizon -- up to one year. The
  delivery quantities may vary markedly day to day and season to season,
  and the the underlying routing problem may have rich constraints
  -- e.g., time windows, compatibility constraints, and compartments.

  While many approaches to solve the fleet size and mix (FSM)
  problem in various contexts have been described, they usually only
  consider a ``representative day'' of demand. We present a method
  that allows any such single-day FSM solver to be used to find
  efficient fleets for a long time horizon.

  We propose a heuristic based on column generation. The method
  chooses a fleet of vehicles to minimise a combination of
  acquisition and running costs over the whole time horizon. This
  leads to fleet choices with much greater fleet utilisation than
  methods that consider only a subset of the available demand data.
  Issues with using a heuristic sub-problem solver within the context
  of column generation are discussed.

  The method is compared to two alternative methods on real-world
  data, and has shown to be very effective.
\end{abstract}


\maketitle

\section{Introduction}
\label{section: introcution}
Efficient transportation is a fundamental component of a nation's
economy, accounting for a significant fraction of the total cost of
products. In order to stay competitive, particularly in the context of
a global and connected market, companies are pressured to reduce their
transportation expenses by making smarter use of their
resources, e.g., fleets and warehouses. Of foremost importance, in
this regard, are the decisions made at the strategic (long-term) and
tactical (mid-term) levels, as they determine how the day-to-day
operations can be carried out, which in turn reflects on the
efficiency of the system. We look at one such tactical-level problem.

Consider the following situation.  A freight company is bidding for a
large contract with a national supermarket chain. If successful, they
will devote a certain number of existing and newly-purchased vehicles
to the new contract. They have one year's worth of demand data from
the supermarket chain, exhibiting large variations in quantity on day to
day and season to season basis. In order to submit a competitive but still
financially viable bid, the freight company must determine the number
and type of vehicles required to most effectively carry out the
deliveries across the whole year. This is the problem we address here.

\emph{Fleet design} is the problem of determining
the size and composition of a fleet of vehicles to carry out the daily
delivery operations of a company. The problem can be informally stated
as follows: given
\begin{enumerate}
\item the demand of a set of customers over a period of
  time, \label{item:demand}
\item a ``catalogue'' of truck types with different characteristics, e.g., capacity,
  running costs, compartments, etc., and \label{item:veh_types}
\item a model describing the constraints and costs incurred in running the
  daily delivery operations, \label{item:VRP}
\end{enumerate}
identify the most efficient fleet to satisfy the demand across the
whole time period.

In our context, the demand (\ref{item:demand}) can be represented
either by historical data \citepar{kilby2016fleet} or forecast data.
The underlying VRP model (\ref{item:VRP}) can be a rich vehicle
routing problem (see 
\cite{caceres2015rich,drexl2012rich,hasle2007industrial}) -- a model
capturing real-world operational constraints such as time windows for
delivery, travelling distances, capacities, compatibilities, driver
breaks, and so forth.

Most of the existing literature considers fleet design as an extension
of the classic vehicle routing problem where the fleet composition is
not an input to the problem, but rather a decision to be
made. There are several variants of this extension which are
usually called Fleet Size and Mix VRP or Heterogeneous VRP (according
to \cite{baldacci2008routing}) depending on whether the fleet is
unlimited or not. We will use the acronym FSM to denote this class of
problem based on the classical VRP.
Due to its huge search space, the FSM is much harder to solve than its
VRP counterpart. For this reason, the FSMs models rarely extend beyond
considering a single representative day of demand data, or model
not-so-rich VRP variants.

The main features that define the VRP variant studied here are
\begin{itemize}
\item day-by-day \footnote{
  Note that we use the term \emph{day} because in most
  contexts, the VRPs are solved on a per-day basis. However, the term
  should be understood to mean \emph{a suitable unit of time} in the
  context of the problem at hand.}
  demand known over a long time horizon;
\item potentially, a rich underlying VRP; and
\item we wish to determine fleet size and mix.
\end{itemize}
We call this problem the Long-Horizon Fleet Size and Mix (LHFSM)
vehicle routing problem to differentiate it from the standard FSM.

While a solver for the (single-day) FSM can often be used to solve a
multi-day instance by careful use of time window constraints, such
approaches do not usually scale well. Time horizons of one year for
hundreds of customers are usually well beyond the scale that can be
handled by these methods.

Another common way to approach the problem is to reduce the problem to a
``representative day'', or even a single ``big day'', and then solve
using that single set of data. However, this approach can not
adequately represent the variability in demand day to day and season
to season. In addition, rich constraints may make a good solution for
one day infeasible for another. For example, a good solution for a
``big day'' scenario may use many large trucks. But vehicle
compatibility constraints, which forbid large vehicles visiting some
customers, may make that solution expensive, or even infeasible, for a
day with more moderate demand.

The method proposed here provides a framework in which a single-day
FSM solver (appropriate for the underlying single-day rich VRP) can be
used as a sub-problem solver to determine an efficient solution for
the longer horizon. We show that the fleets produced are more
effective than those produced considering only ``representative''
days.

Even though some companies build their fleet slowly over a long
horizon, buying one vehicle at the time, there are cases where a
company needs to design an entire fleet. Instances include fleet
down-sizing after a merger, or fleet acquisition in response to gaining
a new, large contract. A company may also want to determine their best
fleet, so as to have a ``target fleet'' in order to guide purchases
over time.

One notable attempt \citepar{kilby2016fleet} has been made to solve the
fleet design problem for a multi-day rich vehicle routing formulation.
The proposed approach is based on a pre-processing step that
identifies a set of Pareto dominant days. The idea is that if a fleet
can cover a ``big-demand'' day then it will also be able to cover all
``smaller'' days, but less so in the general case.  This is possible
in the specific problem presented because a dominance rule,
\bluemark{i.e. a partial order between days, can be established.  While
the approach proved to be viable for the problem at hand, the authors
identified two major issues: scalability and, more importantly,
generality.  The scalability issue arises from the fact that the whole
multi-day problem associated with the Pareto front has to be solved at
once in order to guarantee that, on each day, a subset of the same
fleet of vehicles is used. Given that the single-day problem is
already \textbf{NP}-hard, solving a multi-day variant for a large
Pareto front can easily
become intractable.  A more fundamental
problem is generality -- that is, the problem-specific pre-processing
technique cannot be extended to arbitrary VRP variants. For example
this approach does not work when compatibility constraints are
considered. The fleet for the Pareto days could, and in fact often
will, be infeasible for days with smaller demand but with more
compatibility constraints. The situation is even more complex when
time related constraints, such as time windows, are
considered. Defining a dominance rule for rich VRPs is not trivial, if
at all possible.
An additional inherent problem is the efficiency of the fleet. Since
the fleet produced will be tailored for big days only, the authors have
shown there will often be a high number of idle vehicles per day,
which is not desirable.
However, in Section \ref{section: pareto} we elaborate on how our
method is able to implicitly identify a subset of days that are the
critical ones. This will be a consequence of the mathematical
framework we present.}

\ignore{This is possible in the specific problem presented
because a dominance rule \redmark{,i.e. a partial order between days,} can be established.  While the approach
proved to be viable for the problem at hand, the authors identified
three major issues: scalability, generality, and the need to consider
hiring options. The
scalability issue arises from the fact that the whole multi-day
problem has to be solved at once in order to guarantee that, on each
day, a subset of the same fleet of vehicles is use{}d. Given that the
single-day problem is already \textbf{NP}-hard, solving the multi-day
variant for long planning horizons, e.g., one year, is
intractable. The generality issue arises from the fact that the
problem-specific pre-processing technique cannot be extended to
arbitrary VRP variants. \redmark{ }
The last issue concerns the efficiency of a fleet; the
results in \cite{kilby2016fleet} have shown that, by dropping the
requirement of serving all the demand using owned vehicles, companies
can achieve a much greater vehicle usage.
The hiring options were not, however, explicitly
considered in the presented approach.
}

In this paper, we aim to address the issues of scalability and
generality. We present a general and scalable framework to determine
efficient fleet configurations from historical or forecast demand
data. Notably, we allow the use of a rich FSM solver to solve the
sub-problems in a column generation scheme.  This means that the
suggested approach can be used effectively, regardless of the richness
of the underlying routing problem or the length of the demand
horizon. The technique can be applied
wherever a solver for the \emph{single day} routing problem is
available. The resulting algorithm is a combination of
heuristic and mathematical programming methods and therefore can be
considered a ``matheuristic''
\citepar{boschetti2009matheuristics}. We show the effectiveness of our
approach on a rich vehicle routing variant representing a fuel
distribution problem for which we have real-world remand data.

The paper is organised as follows: in Section \ref{section:
literature} we review the literature related to fleet design. In
Section \ref{section: problem} we formally define the problem that we
tackle in this work. Sections \ref{section: model} and \ref{section:
method} are devoted to the description of the model and the solution
method proposed. In \ref{section: implementation details} we give a
detailed account of the algorithm implementation. Experimental results and
analysis are presented in Section \ref{section: simulations}. We
finally conclude with some remarks and future research directions.

\section{Related Work}
\label{section: literature}

In the following, we briefly review the existing literature on fleet
design and rich vehicle routing.


In \cite{crainic2003long}, the author provides a complete survey of
long-haul freight transportation. The major decision are divided in
three levels: strategic, tactical and operational. Examples of such
decisions are: demand modelling and location of facilities for the
strategic level, service network design and resource acquisition for
the tactical level, routing and loading of vehicles for the
operational level. A key point is that decisions at each different
level have a strong influence on the decisions involved in next level
and provide information for the decision at the previous level.  \\ In
\cite{hoff2010industrial}, the authors present a thorough survey of the
literature related to the problem of managing a fleet, with an
emphasis on
industrial applications. Following their classification the problem can
be looked at from three different points of view, which correspond to
the levels mentioned in \cite{crainic2003long}: strategic, tactical
and operational.  The main difference is the time horizon considered.


We concentrate our survey on papers addressing a
tactical-level problem such as the one studied here. Strategic models
usually focus on the evolution of the fleet through time and do not include
any routing component in the objective, and are therefore addressing quite
different problems. Operational level models (surveyed in
\cite{baldacci2008routing,hoff2010industrial}) do not usually scale to
the lengths of time horizon we are considering, and so are also
omitted.

As noted in \cite{hoff2010industrial} and
\cite{salhi1993incorporating}, the literature has not addressed
tactical level problems in as much depth as the operational level.
We quote \cite{hoff2010industrial}:
\vspace{3mm}
\begin{quoting}
  \textit{
    A large part of the literature focuses on operational questions,
    along the line of ``what to do given a certain fleet mix and a
    given set of service request''. This is in contrast to the more
    tactical, or strategic, ``which vehicles should we acquire to best
    solve our daily routing problem for the next half year''. There is
    a big absence of papers addressing these more tactical questions,
    and also on how to make robust or resilient solutions. }
\end{quoting}
\vspace{3mm}

Two exceptions are:
\citepar{yoshizaki2009scatter}, where the authors consider a week as an
horizon, although eventually solve each day separately, and
\citepar{salhi1993incorporating}, where the authors note the same lack
in the literature and develop an advanced heuristic. Their solution
method is based on a route perturbation procedure. After initially
solving a classical VRP with a given capacity to create a starting
solution, the algorithm checks which vehicles can be used on the routes
and whether it is better to remove a route and insert the customers in
other routes. Other refinement techniques are then applied to improve
the solution. The whole algorithm is repeated several times with other
starting solutions. The overall performance is good and the algorithm
proved to be stable and able to handle complex constraints but is not
tested on long horizons.

\redmark{As we noted before, a first paper attempting to solve on a year worth of data is
\cite{kilby2016fleet}. The authors make use of the particular structure of the specific problem considered to identify a Pareto dominant subset of days. They solve the problem only for the identified set by means of adaptive large neighbourhoods search. While this approach solves a (small) number of days simultaneously, it can still fail to produce a good fleet, for instance in the presence of compatibility constraints.}

\ignore {A first paper attempt to solve on a year worth of data is
\cite{kilby2016fleet}. Here the authors consider a rich VRP and
solve the problem on a Pareto dominant subset of days by means of
adaptive large neighbourhoods search. The particular structure of the
problem allows the authors to guarantee that the so obtained solution
covers all days in the problem. To obtain a more flexible fleet the
authors solve the problem on various subsets of days, corresponding to
layered Pareto fronts. While this approach solves a (small)
number of days simultaneously, it can still fail to produce a good
fleet, for instance in the presence of compatibility constraints.
}

We also differentiate the problem studied here from two classes
of problems involving multi-day horizons: the periodic vehicle routing
problem (PVRP, \cite{francis2008period}) and the inventory routing
problem (IRP, \cite{coelho2013thirty}). In the PVRP, customers are
visited a fixed number of times over the planning period. While similar in
some ways, the LHFSM is significantly different for two
reasons. First, the fleet selection component is not present in
PVRP. Second, in the PVRP the delivery day is a decision variable,
whereas in LHFSM the delivery day is fixed.

In the IRP, we know the rate of demand at each customer, and must
visit customers a number of times over the planning horizon,
delivering sufficient goods to avoid a stock-out. Again in the IRP,
the delivery day is a decision variable, and no fleet selection is
required.  Also in the LHFSM, we do not have a demand rate
available. \bluemark{Note that in both the PRP and IRP, the complexity
lies in the choice of an efficient schedule for each customer and on
how this combines with the routing plan. In the LHFSM, on the other hand,
there is no such problem: the schedule is fixed. What links the days is
only the choice of the fleet.}

Another class of problems involving multi-day scenarios are the ones studied in Stochastic Optimisation
\citepar{Gend14Stochastic}. One method is to sample a number of
scenarios (days), and design a set of routes where the objective is to
minimise the sum of costs across the all scenarios. In the most
general case, all 365 days would be used to evaluate the objective,
although this would be computationally very expensive. Another method
is to analyse the demand patterns of each customer, and fit a
probability distribution function. We can then use this to determine
the probability of needing to visit a customer, and the probability of
exceeding capacity. However, in both cases, this in effect imposes a
``grand tour'' constraint: customers, if they require service, are
visited in the same order every day. Moreover, if time windows are
present, the grand tour must observe the time windows in each scenario
(or in expectation), further constraining the routes. If such a
``grand tour'' constraint is not present in our problem, enforcing it
can lead us to very sub-optimal solutions.

\section{The problem} \label{section: problem}

The problem we address is that of optimising the costs of acquiring
and running a fleet of vehicles over a long planning horizon, e.g.,
one or several months, while guaranteeing the necessary level of
service to satisfy the demand of all customers. This requires taking
into account both \emph{fixed} costs (\emph{Capital expenses} or
\emph{CAPEX}), e.g., acquisition and maintenance expenses, as well as
\emph{operational} costs (\emph{Operational expenses} or \emph{OPEX}),
e.g., fuel and drivers' wages.

More formally, let $I$ be the set of days, or \emph{horizon}, that we
are considering, and $T$ be the set of all available vehicle
\emph{types}. Each vehicle type $t \in T$ has a fixed cost $b_t$,
which may aggregate, for instance, acquisition and maintenance
costs. On each day $i \in I$ we must satisfy the demand of a set of
customers while respecting a rich set of constraints, e.g., time
windows, pickup and delivery, compartments, or vehicle-commodity
compatibilities, etc. We denote the rich vehicle routing problem for
day $i \in I$ as $VRP(i)$, and with $FSM(i)$ the corresponding fleet
size and mix problem. As we will show, our fleet design approach can
be applied to any rich vehicle routing problem,
as long as a suitable solver
for $FSM(\cdot)$ is available. This is a key point. To be precise, we
assume to have a suitable solver (which we will refer to as
FSM-solver), either an exact method or a heuristic,
that can solve the FSM problem for a single day while honouring
upper bounds on the number of vehicles of each
type $t \in T$.

Our goal is to find a fleet that minimises the sum of fixed and
operation costs across the whole horizon. In the next section we
present a mathematical programming model for this problem.


It is obvious that if $\| I \| = 1$ then we have nothing to do, as one
run of the $FSM$-solver would be enough. If $\| I \| > 1$ than the
difficulty lies in finding a good balance in the fleet.


\begin{figure}[h]
  \centering
  \subfloat[][\emph{Day 1. Option 1}.]
  {\includegraphics[width=.45\textwidth]{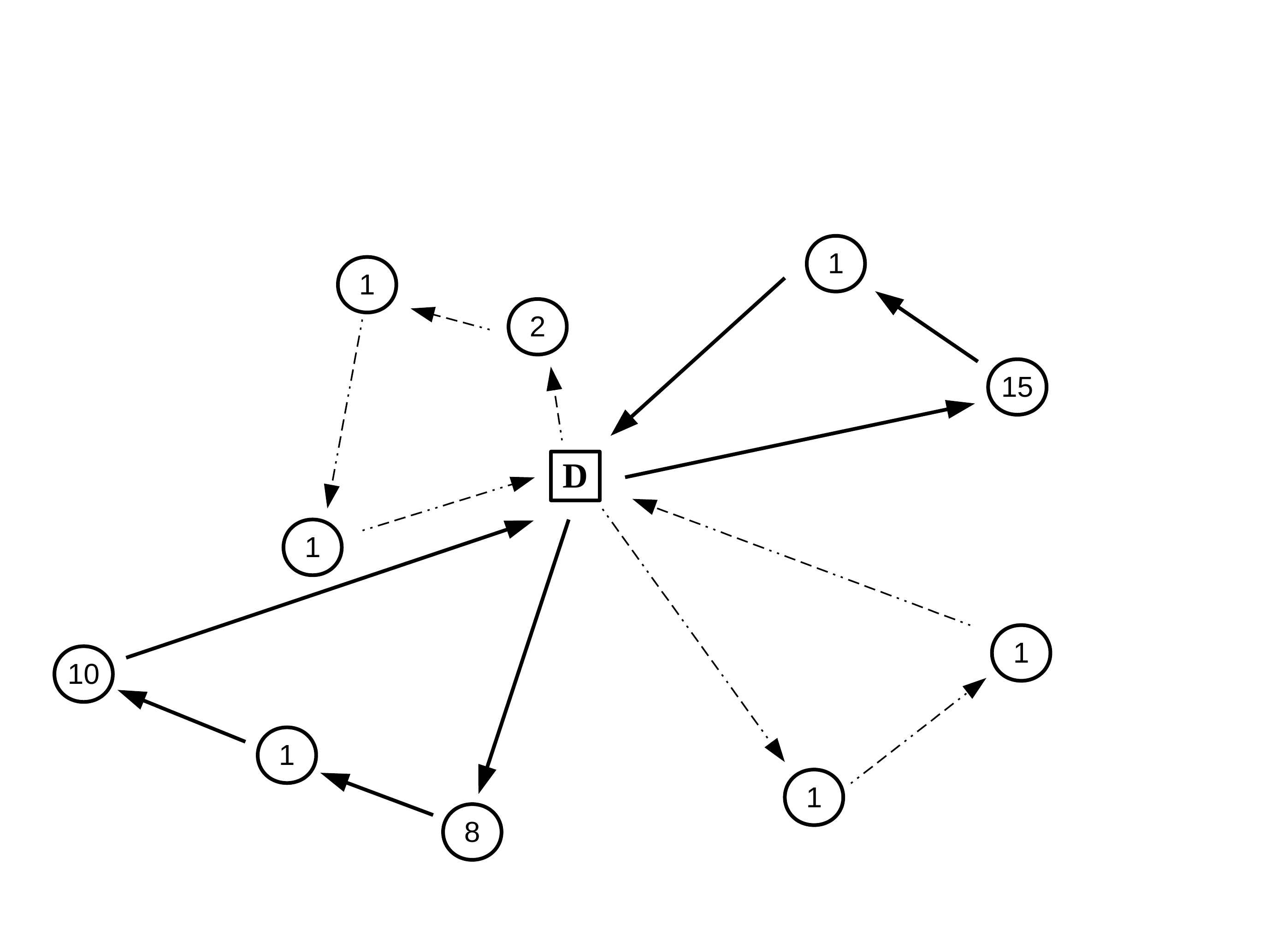}} \quad
  \subfloat[][\emph{Day 1. Option 2}.]
  {\includegraphics[width=.45\textwidth]{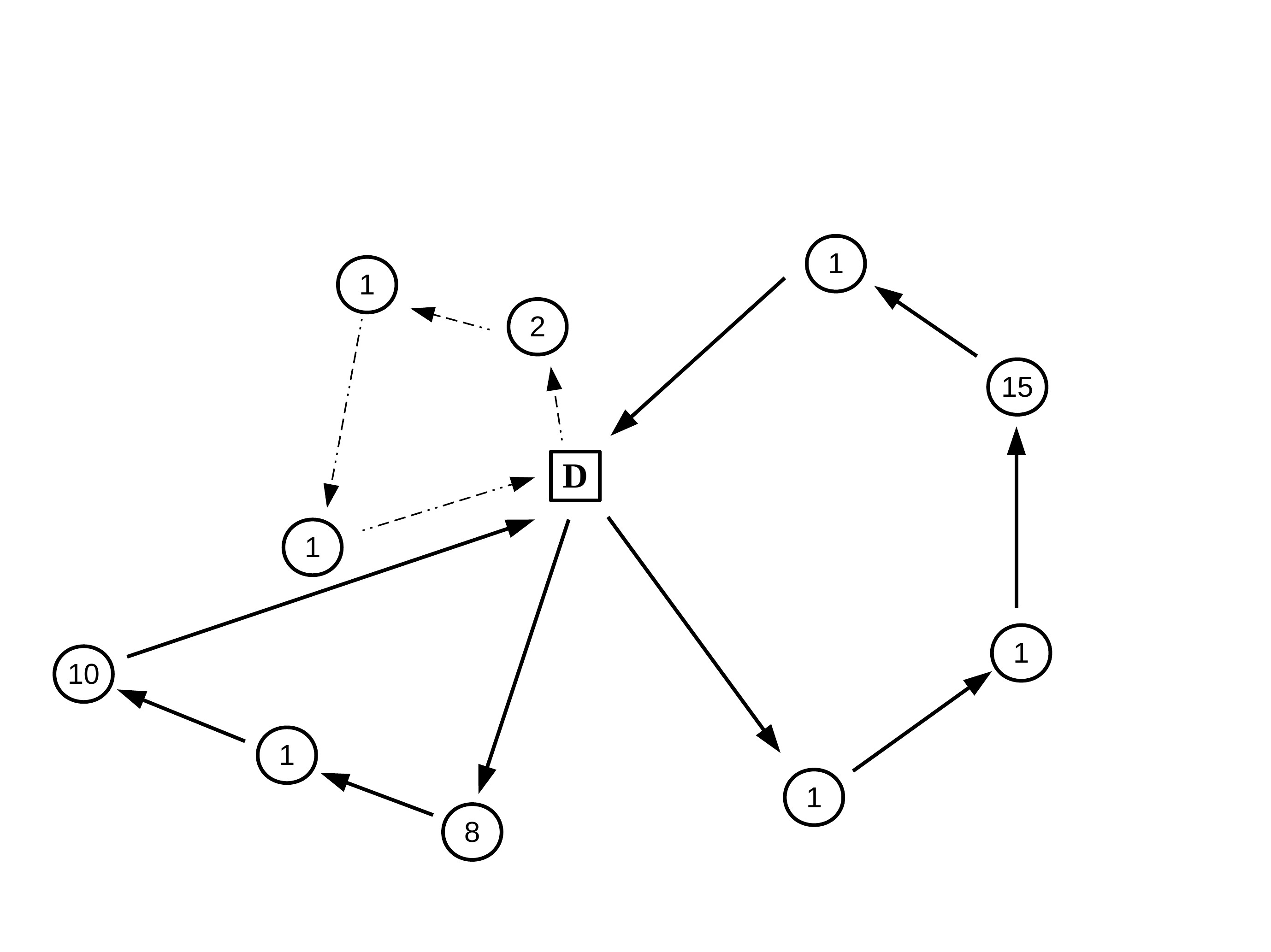}} \\  
  \subfloat[][\emph{Day 2. Option 1}.]
  {\includegraphics[width=.45\textwidth]{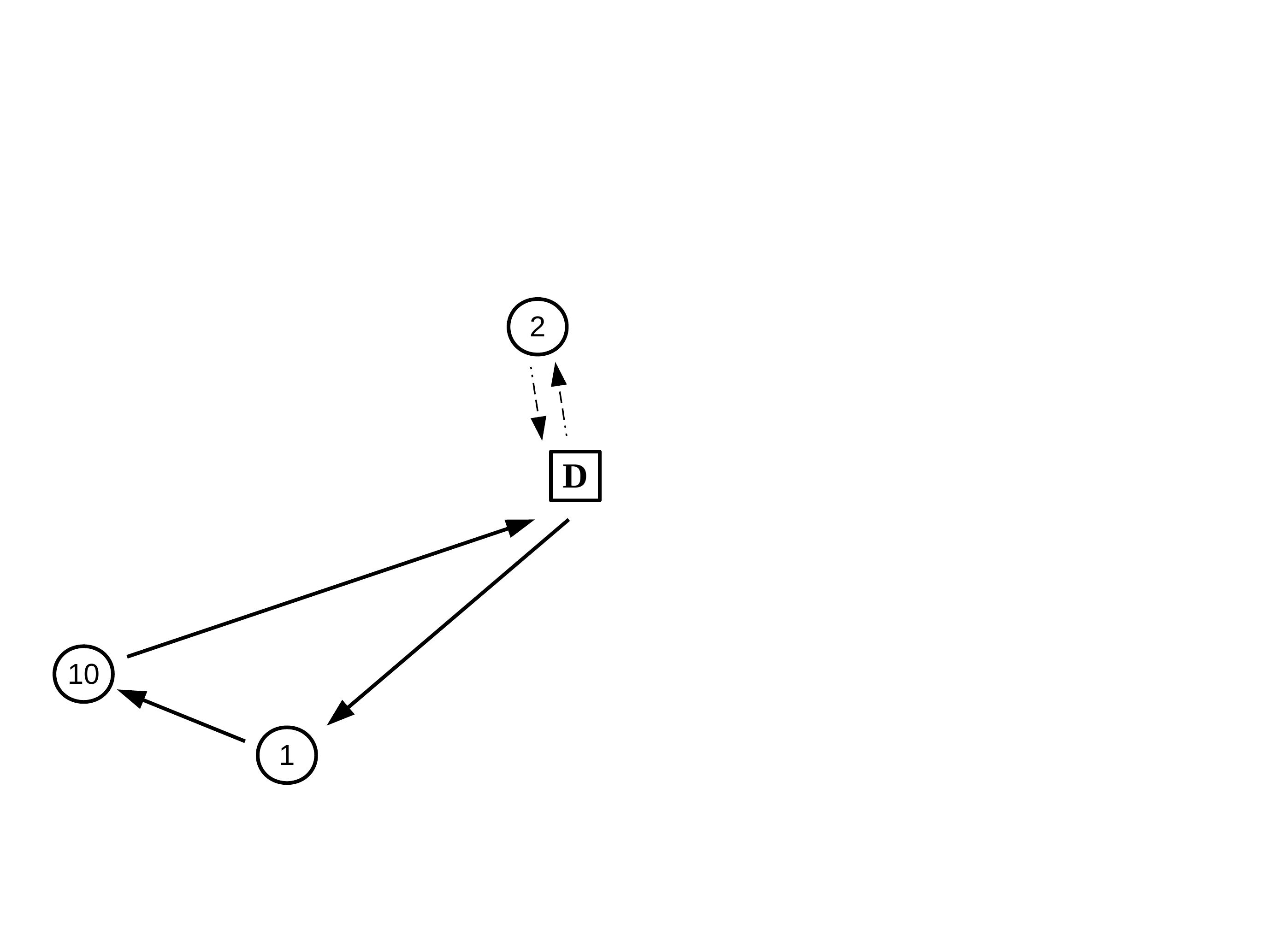}} \quad
  \subfloat[][\emph{Day 2. Option 2}.]
  {\includegraphics[width=.45\textwidth]{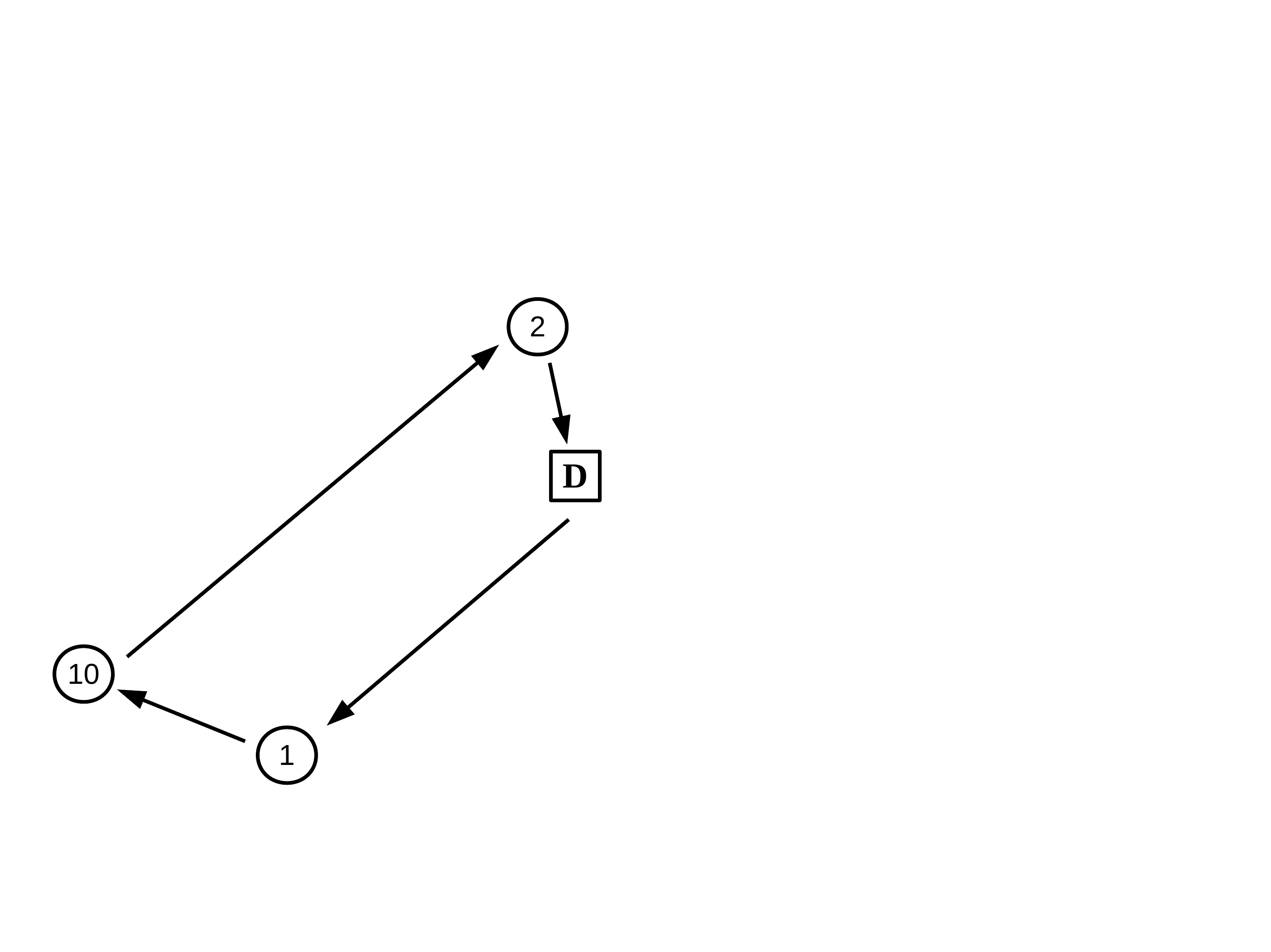}}
  \caption{Suppose
          we have only two types of
          vehicle, with capacity 20 and 5,
          respectively represented by the
          solid and dotted lines. The
          smaller vehicle has a smaller
          routing coefficient. We could
          serve day one with 2 vehicles of
          each type with a routing cost
          of, say, 10. Although we could
          also use 2 vehicle of the first
          type and only 1 of the second,
          which yields a different routing
          cost. Note that option 1 for day
          2 is completely different from
          option 1 for day 1. Moreover the
          fleet used for option 1 of day 2
          would not be feasible for day 1,
          as we would not have enough
          capacity. On the other side the
          fleet used in option 1 for day 1
          is not an option for day 2,
          because we would have 1 vehicle
          left idle.}
  \label{fig: options}
\end{figure}
\setlength{\textfloatsep}{2\baselineskip}

\section{Model} \label{section: model}

\newcommand{\Ni}{N_i}
\newcommand{\rij}{r_{ij}}
\newcommand{\dij}{d_{ij}}
\newcommand{\Fijt}{F_{ij}^t}
\newcommand{\Hijt}{H_{ij}^t}
\newcommand{\Ft}{F_t}
\newcommand{\Et}{E_t}
\newcommand{\Ftp}{\Ft^+}
\newcommand{\Ftm}{\Ft^-}
\newcommand{\Hti}{H_{ti}}
\newcommand{\hti}{h_{ti}}

Given $i \in I$ we consider the set $\Ni$ of possible fleets that can
be used to solve $VRP(i)$. We refer to these as \textit{options} for
$i$. The options for each day will be, in general, different.

Moreover, if we assume that all the vehicles in an option must be used\footnote{\redmark{In other words, there cannot be any idle vehicle in an option.}}, the number of options are guaranteed to be finite.
We denote by $\rij$ the cost of operating option $j \in \Ni$ of day
$i \in I$. This value is calculated excluding any fixed
costs ($b_t$), but including any possible vehicle-specific
coefficients, e.g., fuel cost per kilometre, hourly driver wage, etc.
For every day $i \in I$, every option $j \in \Ni$ can be represented
as a vector of length $|T|$ specifying how many vehicles of each type
are used. We denote by $\Fijt$ the number of vehicles of type $t \in
T$ used by option $j$ for day $i$. Figure~\ref{fig: options} gives 
examples of option sets.

\indent For each option $j \in \Ni$ in each day $i \in I$, we
introduce a binary decision variable $\dij$ representing whether a
particular option is actuated or not. Moreover, the integer decision
variables $F_t$ represent, for each $t \in T$, the number of vehicles
of type $t \in T$ that are part of the overall fleet.\\
\indent We seek to minimise the sum of the fixed ($b_t$) and operation
($\rij$) costs of the fleet across the whole planning horizon. The
model follows.

\begin{align}
  \text{[M] minimise} & \sum_t \Ft b_t + \sum_{i,j} \rij \dij &  \label{eq:m_obj}   \\
  \text{subject to}   & \sum_j \dij = 1                       & \forall \, i \in I                   \label{m_options}                           \\
  & \sum_j \Fijt \dij \leq \Ft            & \forall \, i \in I, t \in T          \label{m_vehicles}      \\
  & \Ft \geq 0, integer                   & \forall \, t \in T                    \label{m_integrality1}  \\
  & \dij \in \{0, 1\}                     & \forall \, i \in I, j \in \Ni \label{m_integrality2}
\end{align}
Constraints (\ref{m_options}) means that we can only select one option
for each day.  Constraints (\ref{m_vehicles}) ensures that the overall
fleet has enough vehicles of type $t$ to actuate the chosen option for
every day. Constraints (\ref{m_integrality1}) and
(\ref{m_integrality2}) enforce the integrality of the solution. The
objective is the sum of fixed and the operational costs.

The solution to this problem will give us an exact solution to LHFSM
if the column costs $\rij$ are calculated exactly. In many situations
where a rich VRP underlies the selection process, $\rij$ will be
calculated by a heuristic method. In this case the above model will
give a heuristic solution to the fleet selection problem.

Note that, provided that $\dij$ is positive, Constraints
\ref{m_options} implies $\dij \leq 1$. Therefore we can replace
Constraints \ref{m_integrality2} with $\dij \geq 0, integer$. We
decided to keep $\dij$ binary to highlight its semantics as binary
choice. Similarly we could replace Constraint \ref{m_options} by
$\sum_j \dij \geq 1$ due to solving a
minimisation problem.

\section{Solution method} \label{section: method}

Every coefficient $\rij$ in (\ref{eq:m_obj}) must be generated by
solving the the corresponding $VRP(i),\, i \in I$ using the fleet $j
\in \Ni$. Each one of these $VRP(\cdot)$ is \textbf{NP}-hard. Moreover
the number of days $|I|$ can be quite large, e.g., one or more
years. Determining the number of possible options $|\Ni|$ would
already be a hard problem, consequently the full enumeration of the
columns is computationally intractable and solving the linear
relaxation of [M] is not a viable option. \redmark{In these cases it is common to use a column generation (CG) approach. However CG can only be used to solve linear problems. In order to solve integer programs it is necessary to embed the CG method in other algorithms aimed at solving the integer version of the problem. We proposed two such approaches later in the section.}



In the following, we first briefly review the column generation
technique in general terms, then we describe the necessary steps to
apply it to our model.

\subsection{Column generation}

Column generation \citepar{lubbecke2005selected} is a technique for
solving linear programs where the number of variables, or equivalently
columns, is intractable or hard to enumerate. The approach exploits
the idea that only a small number of variables are needed to determine
the optimal solution, as most of them will assume value zero.

The first step is to decompose the original problem into a
\emph{master} problem and a \emph{sub-problem}. At each iteration, the
master problem only considers a restricted number of variables, or
equivalently columns, and it is solved to optimality. The sub-problem
is an artificially defined problem that tries to identify which new
columns are to be added to the master problem in the next
iteration. The objective function of the sub-problem is the
\emph{reduced cost} of the new columns with respect to the optimal
dual variables of the current optimal master solution. The constraints
in the sub-problem ensure that a solution is indeed a valid column in
the original problem. This phase is commonly referred to as
\emph{pricing} of columns. If we identify a column with negative
reduced cost, we add it to the master and proceed to the next
iteration, otherwise we stop. 

\makeatletter
\renewcommand{\ALG@name}{Column Generation}
\makeatother
\renewcommand{\thealgorithm}{}

\begin{algorithm*} [ht]
\caption{}
  \begin{algorithmic}[1]
    \State Solve the (restricted) master problem.
    \State Retrieve the dual variables from the master.
    \State Solve the pricing problem to identify a new variable. \label{step: SP}
    \State If the objective of pricing problem is non-negative we have reached optimality. \label{step: optimality}
    \State Otherwise add the identified variable to the master and go to Step 1.
  \end{algorithmic}
\end{algorithm*}

A fundamental aspect in the implementation of a column generation
scheme, is to recognise the structure of the sub-problem, and to
identify a problem-specific technique to solve it efficiently.
For instance, when solving a standard VRP by column generation
\citepar{desrosiers1984routing}, the sub-problem can be seen as a
constrained shortest path problem which can be solved by labelling
algorithm. 
\redmark{Note that it is not necessary to identify a single
sub-problem. A problem can be decomposed in a master and several
sub-problems. In fact this happens for our model [LM]. The
extension of CG to this scenario is immediate. In Step \ref{step: SP}
we simply solve all sub-problems and add several variables to the master. The optimality condition (Step \ref{step: optimality}) translates to all sub-problems having a non-negative objective. }

\subsection{Sub-problem}

In order to apply column generation to our problem we consider [LM], the
linear relaxation of [M], and denote by $p_i$ and $q_{ti}$ the dual
variables of Constraints \ref{m_options} and \ref{m_vehicles},
respectively.  \redmark{It is evident that we have one sub-problem for
each day $i \in I$. The sub-problem associated with day $i$ aims at
generating a variable $\dij$ with negative reduced cost. We now
illustrate the special structure of the sub-problems.}

The reduced cost of a variable $\dij$ is
\begin{equation*}
  \rij + \sum_t \Fijt q_{ti} - p_i.
\end{equation*}

Each column is uniquely identified by $\rij$ and the $\Fijt$s,
which represent variables of the sub-problem. Ignoring the term $p_i$,
which is a constant in this context, we observe that the reduced cost
is exactly the objective of an FSM where the fixed costs of the
vehicles are represented by the $q_{ti}$ variables and $\rij$ represents the
cost of optimally operating the fleet determined by the $\Fijt$
variables. \redmark{In other words, sub-problem associated with day $i$ generates a new option $j \in N_i$.}

Since the FSM is \textbf{NP}-hard, optimally solving the sub-problem
efficiently for large problems is usually not viable. Instead, we
allow the use of heuristic methods which are able to obtain good
quality solutions in a short time. Using a heuristic method
to solve the sub-problem makes our algorithm a heuristic. However, the
method described still converges to a good solution to the original
problem, modulated only by the quality of the underlying FSM
heuristic.  We regard this as one of the main strengths of our
approach. If a solver for the single-day version of the problem is
available, we can use it to solve the multi-day variant of the same
problem.


Several issues arise when using a heuristic solver. The primary one is
that, due to the stochastic nature of many heuristic solvers,
different solutions can be produced for the same input. Thus, an
$\rij$ for given a fleet $j$ on day $i$ may be calculated at one
time. Then, later in the execution, a new, better $\rij'$ may be
produced. When this happens, it is enough to replace the original
$\rij$ with the new $\rij'$ in the master problem. Because we only
ever reduce the costs of columns, this does not effect the convergence
of the column generation.

\subsection{\redmark{Proposed Methods}}
\label{section: algorithms}
As previously noted, CG only solves a linear problem. In order to solve [M] we propose two approaches. In the first one we use CG to solve [LM], and then we solve a restricted version of [M] involving only the columns generated in the linear phase. This method is called \textit{restricted master heuristic} \citepar{joncour2010column}. The second approach consists in a branch \& bound tree search where CG is used at each node to provide lower bounds. This type of method is generally called \emph{branch \& price}. We denote the two methods by, respectively, RMH and BAP.

Note that RMH is heuristic in its nature, as the set of column that identify an optimal basis for the linear problem [LM] might not contain the optimal set of columns for the integer problem [M]. Oppositely, \emph{branch \& price} approaches aim to be exact. However, as we solve the sub-problem by means of an heuristic, both the methods are heuristics and do not guarantee optimality.

\subsection{Historical and Stochastic Data}

The demand over the time horizon of each customer is
considered to be known and given as an input to the
problem. However, the approach described can also be used
effectively to handle stochastic demand for a single-day problem.
Given some expression of the stochastic demand at each customer,
an artificial multi-day problem could be created by generating
different realisations of that demand to produce multiple
scenarios. The method described above can then be used to determine
a fleet that can cover the variable demand most effectively.
While this approach seems to have some strength, we have not
investigated it further.

\subsection{\redmark{An Extension of the Pareto Approach}}
\label{section: pareto}
In the Introduction we commented on the difficulty of extending the method proposed in \cite{kilby2016fleet} to rich vehicle routing problems. A byproduct of the mathematical programming framework that we are using is that we have some insights on how to extend the idea of identifying a subset of critical days. 
Let us define the dual of problem [LM].

\begin{align}
  \text{[D] maximise} &   \sum_i p_i                                          &  \nonumber \\
  \text{subject to}   &   \sum_i q_{ti}  \leq b_t  \quad                      && \forall \,  t \in T \nonumber \\
                      &   p_i - \sum_t F_{ij}^t q_{ti} \leq r_{ij} \quad      && \forall \, i \in I, j \in \Ni \nonumber \\
                      &   p_i \in \mathbb{R} \quad                            && \forall \, i \in I \nonumber \\
                      &   q_{ti} \geq 0 \quad                                 && \forall \, i \in I, t \in T\,  \nonumber 
\end{align}

The dual problem ``distributes'' the costs $b_t$ over the days (if $q_{ti}$ are interpreted as daily fixed costs) with
the goal of maximising the sum of the lower bounds ($p_i$ in this context) on each day's total
cost. If, for a fixed $i \in I$ and $t \in T$, Constraint \ref{m_vehicles} is not tight, the Complementary Slackness Theorem implies that $q_{ti} = 0$ and therefore the vehicle is ``free'' in the sub-problem. 
Intuitively, the variables $q_{ti}$ tell us how much a vehicle
of type $t$ is worth on day $i$ if we can only use the options
(columns in the master) provided. In fact, in each iteration of CG, there will be a few days, having some of the $q_{ti}$ that are non zero. This information point us to the days that need other options in order to reduce the overall fixed cost and it is of fundamental importance for an efficient implementation as is discussed in next section. In some sense the dual variables identify a subset of critical days. Notably, this set changes at each iterations, as the insertion of new variables in the master changes the dual problem [D]. However, the main advantage of our decomposition is that we can solve each day independently. Oppositely, in \cite{kilby2016fleet} each Pareto subset has to be solved as a whole, considerably increasing the complexity of the method.


\section{Implementation Details}
\label{section: implementation details}

We now go through the details of our implementation explaining how we
address standard issues that arise when implementing column generation
approaches and new ones originating from the fact that we are not
solving the sub-problems to optimality.
\vspace{3mm}

\paragraph{Branching strategy.}
In a branch \& price scheme one has to carefully choose the branching
strategy so that the sub-problem structure is preserved. Once the
linear relaxation of the root node is solved, we choose a fractional variable
$F_t$ and branch on it by imposing, for
example, $F_t \leq m_t$ and $F_t \geq M_t$ where $m_t$ and $M_t$ are
defined as $\lfloor F_t \rfloor$ and $\lfloor F_t \rfloor + 1$.
Note that adding
this type of constraint to the master problem does not affect the
problem structure, and constitutes additional information that can be
passed to the underlying FSM-solver. Indeed bounding the maximum
number of vehicles of type $t \in T$ is standard in FSMs, while
limiting it from below can be done by considering $M_t +1$ free
vehicles of type $t$. Consequently a node is identified by the lower
and upper bounds imposed on each vehicle type.

When exploring the branch and bound tree, we choose the node with the
lowest integer objective, and which forbids the currently best-known
integer solution. This has the effect of forcing the algorithm, and
the underlying FSM-solver, to look at different parts of the search
space. Note that, provided we have one
option per day, problem [M] always has a feasible solution.

\vspace{3mm}
\paragraph{Theoretical caveats.}
As pointed out earlier, the fact that the sub-problem is solved
heuristically creates some theoretical issues that, however, can be
easily addressed in the implementation. The first thing we note is
that the heuristic might find the same column (fleet option) multiple
times with different routing costs. If the cost is lower than the
previously found one, we just replace the existing column with the
newly found one. Otherwise we discard it. The main theoretical
inconsistency comes from the fact that the solution of the master
problem at the root node is not guaranteed to be optimal and therefore
not guaranteed to be a lower bound.
Therefore, at some point in the branching tree we might find
a node which has a better linear (and even integer) objective than the
root node's linear objective.
This can happen when the branching constraints force the underlying FSM-solver to look
in parts of the search space that were previously missed. If a better
bound is found during execution, it does not affect convergence, only
execution time. In
effect, nodes that might have been pruned earlier will be pruned when
the new, improved solution is found.

\vspace{3mm}
\paragraph{Initialisation.}
As noted earlier, a good feature of model [M] is that, provided that
we have at least one option for each day, the problem is always
feasible. We therefore need to initialise the problem by finding a
first option for each day. Our choice was to use the underlying
FSM-solver to find the one that generates the best routing option. In
practice, we solve each day assuming that all the vehicle types are
unlimited and that the acquisition cost is zero. We know that such
option will have the best routing cost possible for that day. This
will turn out to be useful during the sub-problem selection and the
lower bound computation, as explained later in the section.

\vspace{3mm}
\paragraph{Number of solved sub-problems at each step.}
One particular feature of our problem is the high number of
sub-problems that have to be solved. As we want to apply the method to
long term scenarios, the number of days (sub-problems) can be in the
hundreds.

It is clear that, in our case, solving the sub-problems is the most
computationally expensive part of the algorithm. Simulations show that
it amounts to $99.9\%$ of the total run time. On the contrary, the
solution of the master problem, integer or linear, is almost
instantaneous. Therefore we are not concerned with the number of
columns we insert in the master problem.

It is therefore advantageous to choose a subset of sub-problems to
solve at each iteration. A similar problem appears in
\cite{gamache1999column} where the number of sub-problems is $\sim
300$. The authors' solution is to order the sub-problems and solve only
a fraction of them ($15$ to $20$) following a cyclic order until one
column with negative reduced cost is found. In our problem not all the
days have the same impact on the final solution; one of the
interesting aspects of the problem lies in the fact that many days are
different from each other. Therefore we will have some days that are
easy to accommodate with a different fleet and others that require more
work.

\redmark{As noted in Section \ref{section: pareto} the dual variables play a key role in the choice of sub-problems to be solved.}
We observe that, if $q_{ti} = 0,\, \forall \, t \in T$ for a given day $i \in I$, then we do not need
to solve the respective sub-problem, since the solution found in the initialisation phase is
already the best and its reduced cost is zero. 
A standard approach would be to compute a lower bound for each sub-problem to estimate the objective. Despite the intuitive appeal this is not a viable option in our problem. Indeed, the sub-problem is a FSM problem, possibly based on a rich VRP. Computing a good lower bound for this type of problem is not trivial, nor computationally efficient.

Instead, we aim to diversify the search, by avoiding solving a sub-problem similar to
ones we have solved before. The measure we adopt is to compare the
dual variables $q_{ti}$  with those used last time we solved the
sub-problem on day $i$, and choose the sub-problem with the greatest
difference. 

To rank the days we therefore use a weighted
total variation. We fix a day $i \in I$ and consider all the options
$j \in \Ni$ that we have found so far. Denote by $q_{ti}^j$ the dual
variables that were used when we found option $j$. Also denote by
$|F_{ij}|$ the total number of vehicles in option $j$ and by
$\widetilde{N}_i$ the number of options found so far for day $i$. The
score of a day is then given by the following formula
\begin{align*}
  \sum_j \sum_t (q_{ti} - q_{ti}^j)^2 \frac{ F_{ij}^t}{\|F_{ij}\|} \frac{1}{\widetilde{N}_i}
\end{align*}

We take the total variation with respect to the variables $q_{ti}^j$
weighted by the importance that vehicle type $t$ had in the found
option. In the current implementation the number of sub-problems that
we solve simultaneously at each iteration is equal
to the degree of parallelism (independent cores or virtual cores) of
the hardware onto which we run the algorithm. Note that, in general, this does
not have to be the case, and that the choice of this number influences more
than just the run time of the algorithm.

\redmark{Alternatively, one could define the score as the minimum variation (i.e., replacing the sum over all options with a minimum) or simply consider the sub-problems whose dual prices are the highest. However, experiments showed that taking the total variation produced better results.}

\vspace{3mm}
\paragraph{Degeneracy and tailing-off effect.}
In cases where the primal problem exhibits high degeneracy
and many equivalent solutions, solving the dual problem provides
better numerical stability, and allows the direct calculation of dual
variables. Therefore, at each master iteration, we directly solve the (restricted) problem [D].

Another recurring problem in column generation is the tailing-off effect. In the last
iterations the columns found have negative reduced costs close to
zero, causing small or possibly null improvement in the
objective. This is due to the fact that different fleets can produce
almost identical routing costs, and hence the solutions we find in the
final stages only decrease the routing cost by a very small
amount. We use early termination, as described in the next paragraph, to stop the tailing off
effect.

\vspace{3mm}
\paragraph{Early Termination Strategy.}
A common approach to early terminate and possibly fathom nodes is to
compute lower bounds of the nodes' optimal objective. The standard
technique in to use Lagrangian lower bounds. Given a node in the
search tree, we denote by $z$ and $\overline{z}$ the optimal value of
the linear and integer master problem, respectively. Similarly, we
denote by $z_{RMP}$ and $\overline{z}_{RMP}$ the optimal objectives of
the current restricted master problem. Finally, we denote by $rc^*_i$
the reduced cost of the exact solution of sub-problem $i$ obtained
using the dual variables of the current CG iteration. The Lagrangian
lower bound (see \cite{lubbecke2005selected}), is

\begin{equation*}
  z_{RMP} + \sum_i rc^*_i \leq z \leq \overline{z} \leq \overline{z}_{RMP}
\end{equation*}

However, as previously noted, we only solve a subset of sub-problems
at each iterations, hence we do not know $rc^*_i$ for all $i \in
I$. Moreover, due to the fact that we are solving the sub-problem by
means of an heuristic, we have no guarantee we found the optimal
solution of the sub-problems, i.e. we only have an estimate of the
lowest possible reduced costs. We propose a simple modification of the
above bounds to overcome these two limitation, and we make use of
these only to terminate a node's execution early. This does not affect
the convergence of the algorithm, but only the way we
explore the search tree.

We explained earlier that those days with $q_{ti} = 0 ,\, \forall \, t \in T$
necessarily have zero reduced cost. We therefore need to further lower bound $rc^*_i$ for all days associated with sub-problems that have at least a positive dual variable $q_{ti}$ and that we did not solve in the current iteration. We recall that

\begin{equation*}
  rc_i^* = r_{ij} + \sum_t q_{ti}F_{ij}^t - p_i.
\end{equation*}


Obviously a good lower bound has to be dependent on the specific underlying VRP variant. The technique we propose is however easily adjustable. For now, let us suppose we only have one commodity, capacity and some other routing constraints. Note that, given a day $i$, any option $j \in N_i$ satisfies the constraint

\begin{equation*}
  \sum_t c_t F_{ij}^t \geq TD_i
\end{equation*}
where $c_t$ is the capacity of vehicle type $t$ and $TD_i$ is the
total demand of day $i$. Let us assume, for now, that $r_{i0} \leq r_{ij}$. We have

\begin{equation*}
  rc_i^* \geq r_{i0} + \text{min} \left\{ \sum_t q_{ti}F_t : \sum_t c_t F_{ij}^t \geq TD_i \right\} - p_i = rc_i.
\end{equation*}

The variable in the problems between the brackets are the integer and
non-negative $ F_{ij}^t, \, \forall t \in T$. The numbers $rc_i$ are
very easy to obtain. Note that if the underlying VRP has more than one
commodity the lower bound computation \redmark{can be improved by adding a constraint for each commodity}. This is also
the case for other constraints, such as compartments or
compatibilities between products and vehicles. Although the same idea
can be easily adapted to such constraints. Let us assume that the
$rc_i$'s are obtained through a similar valid technique. Now denote by
$S$ the set of days $i$ that we solved in the pricing phase and set
$rc_i = 0$ if $i$ is such that $q_{ti} = 0 \, \forall \, t \in T$. We have

\begin{equation*}
  z_{RMP} + \sum_{i \in S} rc^*_i + \sum_{i \in S^C} rc_i  \leq z \leq \overline{z}.
\end{equation*}

\redmark{To obtain the approximate lower bound we assumed $r_{i0}
\leq r_{ij}$ for all days $i$ and options $j \in N_i$. As discussed
in the \textit{Theoretical Caveats} paragraph, due to the fact the we
are using a heuristic to solve the sub-problems we might find an option
with better routing cost later in the algorithm.  However, since the
lower bounds are only used to terminate early, and not to fathom a
node, the approximation is still admissible.}


\vspace{3mm}
\paragraph{Quick column generation and nested column generation.}
It is standard in column generation to keep a pool of columns formed
by the columns generated in previous iterations that did not
have a negative reduced cost at the time. At every iteration one can check the
reduced cost of the columns in the pool to see if they now have
negative reduced cost. We apply a slightly different approach.
We form a pool consisting of \emph{all} the routes that are
generated for each day. At each iteration, we can formulate
a set partitioning model to solve FSM(i) for each day $i \in I$. The
set partitioning model \citepar{desrosiers1984routing} is the standard
and most effective technique to solve VRP and it has been applied to
FSM too \v{choi2007column}. For all $t\in T, i \in I$ we denote
by $R^i_t$ the set of routes found so far for vehicles of type $t$ on
day $i$. Furthermore, we denote by $C_i$ the set of customers on day
$i$. For each $r \in R_t^i, c \in C_i$ let $c_r$ be the route operational cost,
and let $a_{rc}$ be a binary constant indicating
whether customer $c$ is on route $r$. The models are

\begin{align*}
  \text{[M$_i$] minimise} & \sum_{t\in T}\sum_{r\in R^i_t} (c_r + q_{ti})\, x_r  & \\
  \text{subject to}   & \sum_r a_{rc} x_r = 1                       & \forall \, c \in C_i           \nonumber          \\
  & x_r \in \{0, 1\}                     & \forall \, r \in R^t_i, \forall t \in T  \nonumber
\end{align*}

Given that the number of routes does not grow too fast, we can solve
these models to optimality very quickly. If the objective is smaller
than $p_i$ we have found a new column. Otherwise we have to run the
normal column generation process described above.

This technique allows us to rearrange the routes we already have
to generate a new fleet size option. This is very close to what is
knows as Nested Column Generation (NCG, \cite{mason1998nested}),
except that we do not generate any new columns. It is possible to
use NCG by implementing column generation to solve the
sub-problems. However it is not always easy, \redmark{and for some variants it is not even possible}, to solve a FSM by column
generation.

\section{Experimental Analysis}
\label{section: simulations}
In this section we present some experiments to analyse the performance
of the
proposed method. We have at our disposal some real data provided by a
business partner. The context is a fuel distribution problem in north
Australia. We have data coming from two different distribution centres
(DC), therefore representing two instances of the same problem.
One of these two set of data was used in
\cite{urli2017constraint} to validate their heuristic. The authors
note that the set of instances are rather small. Therefore we chose
to use only the second set of data which is computationally more
challenging. Before presenting the simulation we briefly describe the
underlying VRP and the FSM-solver we used for the single day problem.

\vspace{3mm}
\paragraph{Routing problem and solver}
The problem can be summarised as follow: every day we have a set of
requests to be satisfied by a fleet starting and returning to a single
depot. Each requests is characterised by: the demand, which involves
two different commodities, a time window (morning, afternoon or the
whole day), a service time which depends on the quantity being
delivered, variable for each customer and, possibly, compatibilities
with vehicles types. Each customer can be visited more than once, and
hence this is a ``split delivery'' problem. We
have seven different vehicle types available. Each type is
characterised by: a different number of compartments for which a
maximum capacity for each fuel type is given (each fuel type can have
a different density and a compartment can be filled by only one
product at a time), a fixed annual cost, a variable operation cost,
comprising a per-unit time cost (representing, e.g. salary) and a
per-unit distance cost (representing, e.g., fuel cost). Each vehicle
can perform one or more trips which are constrained to be in a certain
time window. During a working day, a vehicle can return to the depot to
refill (i.e., the problem also features vehicle re-use).
A better analysis of the performance of the FSM-solver on the
single day instances and a detailed description of the problem can be
found in \cite{urli2017constraint}. The authors solve it by means of
a Large Neighbourhood Search \citepar{shaw1998using,ropke2006adaptive}
heuristic where the ``repair'' phase is implemented through
Constraint Programming \citepar{rossi2006handbook} Their heuristic was
developed to support fleet size and mix and is the FSM-solver we use
to solve the sub-problem. We do not go into details of the FSM-solver
implementation; the interested reader can refer to
\cite{urli2017constraint}.

\begin{figure}[h]
  \centering
  {\includegraphics[scale = 0.25]{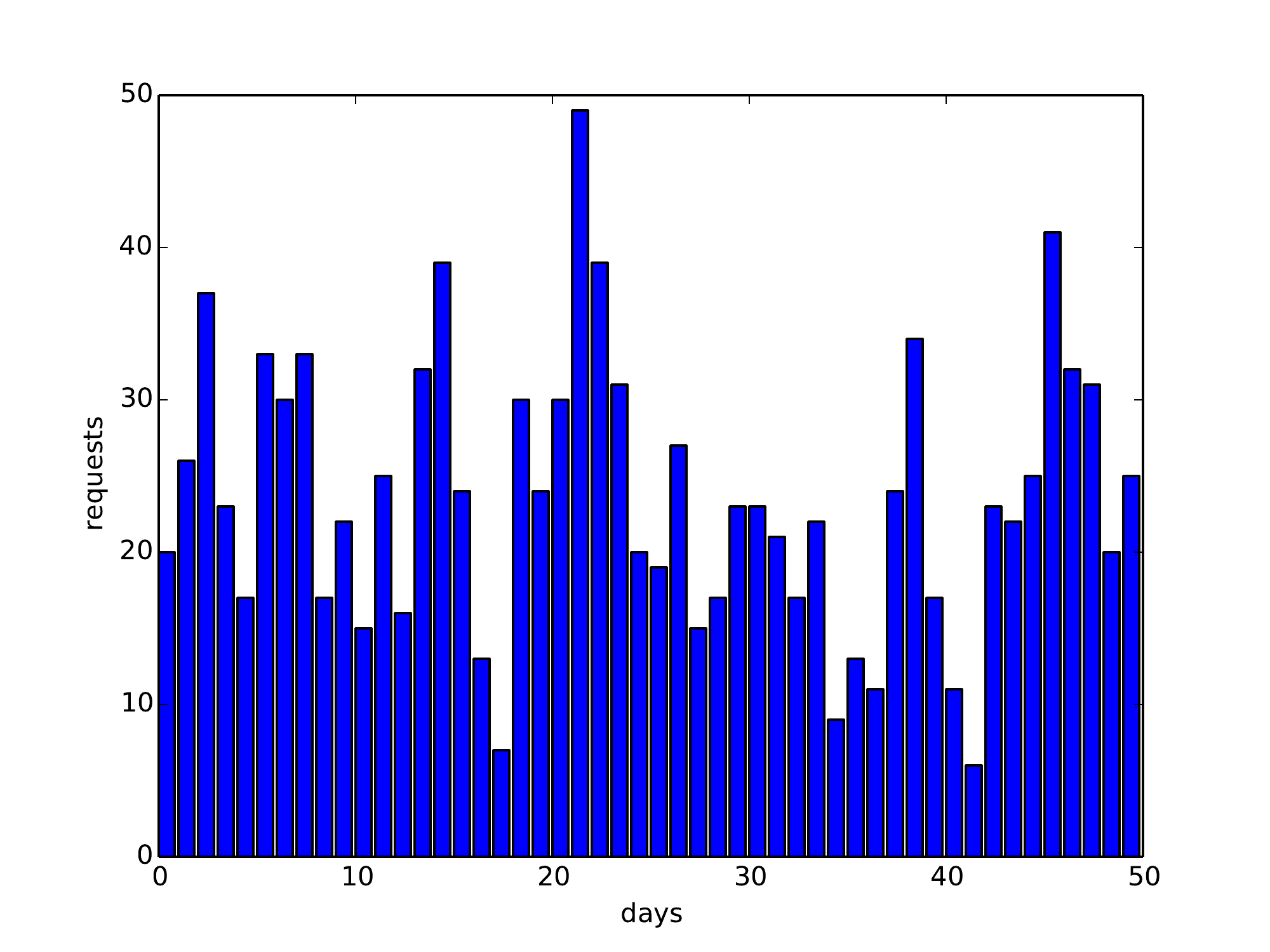}} \quad
  {\includegraphics[scale = 0.25]{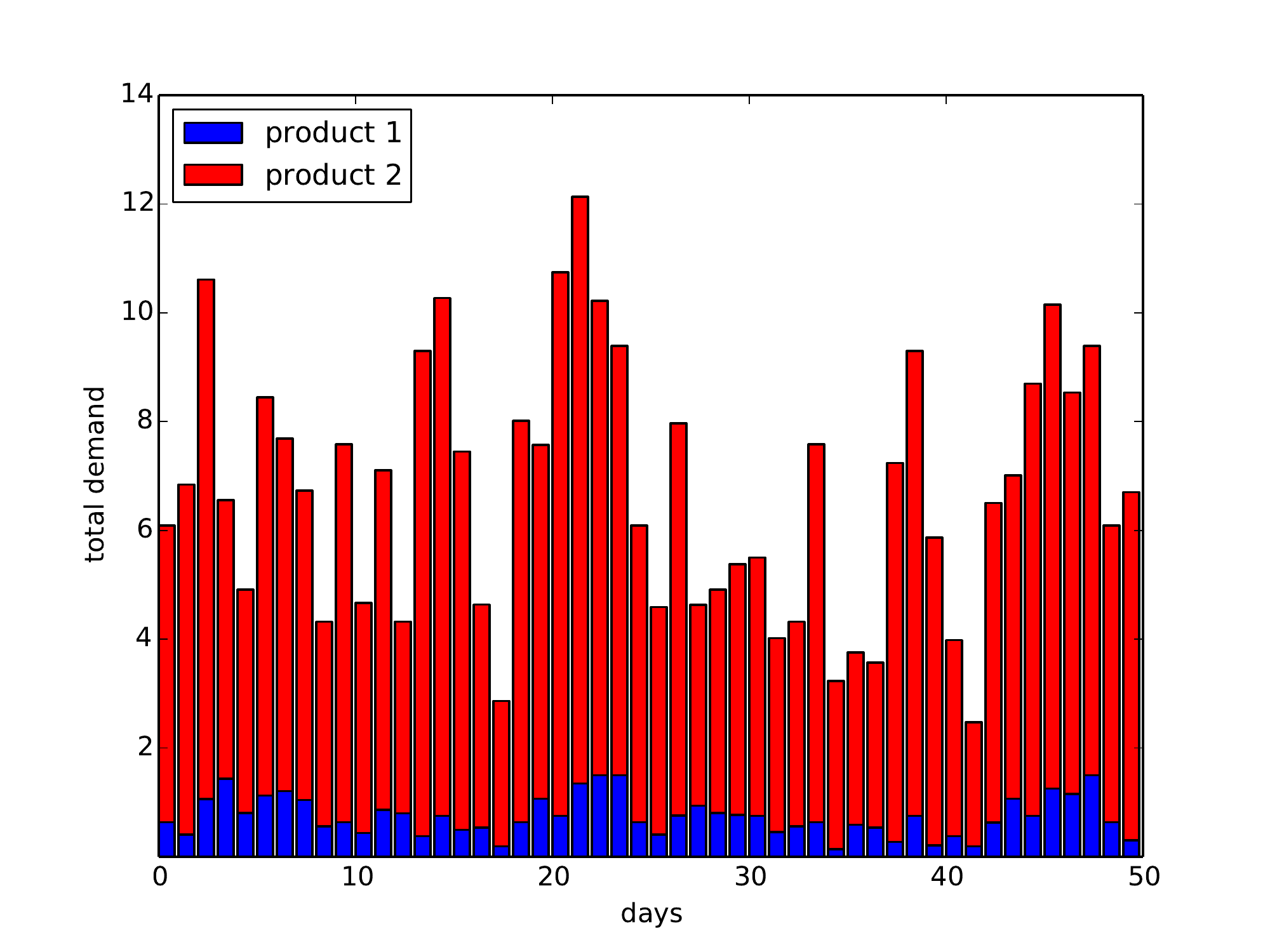}}
  \caption{On the left we plot the number of requests for each day. On the right the total demand for each day for the commodities (the \textit{y}-scale is in hundreds of thousands). The plot includes all the $50$ days.}
  \label{fig: Gladstone data}
\end{figure}
\vspace{3mm}
\paragraph{Available Data}
We have at our disposal a month of demand data. Given that some days do
not have any requests, the actual number of days we have available is
$25$. The different types of vehicle available is $7$ and the number
of customers per day varies in the interval $[20,60]$. To better
analyse the algorithm we modified the data to create a longer
horizon. Since a lot of customers require service on several days it
is possible to perturb the history to create more days. Moreover we
performed a problem reduction by merging a few customers so that the
daily problem is easier and quicker to solve. This resulted in $50$
days with number of customers varying in the interval $[17,52]$. This
gives us the possibility of considering two instances of $25$ days
each. We will denote by DC$1$ and DC$2$ the first and second set of
$25$ days. In order to not compromise the data we tried to maintain
the variation of total demand and number of customers per day. Figure
\ref{fig: Gladstone data} gives a visual representation of how the
daily scenario varies across the planning horizon.

\subsection{Algorithm performance}
\label{sec: algorithm performance}


\redmark{Ideally, we would be able to compare the proposed approaches
with other methods from the literature. We used the method
proposed in \cite{kilby2016fleet} to produce a fleet based on the
days in the first ``Pareto front''. However, due to incompatibility
constraints, this fleet was infeasible for many days in the test
data. Since \cite{kilby2016fleet} does not offer a method to overcome this limitation, we are unable to compare its performance.}


Instead we use an approximate lower bound,
described later, as the ``benchmark'' solution. 

We also compare the performance of the two proposed
methods, RMH and BAP, to other two standard approaches described
in \cite{urli2017constraint}. The first, naive, approach consists in
solving each day separately with amortised fixed costs, i.e.,
$b_t/|I|$ (the overall fixed costs of the vehicles divided by the
number of working days), then taking the fleet that is effectively the
union of the fleet of vehicles obtained by solving each day separately.
The authors name it ``Union
fleet'' (UF) method. We remark that the UF approach treats
every day as a separate entity and therefore is clearly
sub-optimal. However this is the approach that is used in the
literature to solve multi-day FSM problems (see for example
\cite{yoshizaki2009scatter,kopfer2009combining}).

The second approach is to formulate a ``big day'' scenario. The logic
here, often employed in practice, is to form a fleet that can meet the
demand on the day with largest demand, and hence will be able to meet
the demand on all other days.
The method, which we name Subset Algorithm (SA), is illustrated in Algorithm \ref{algorithm: SA}. This is based on the fact that the FSM-solver described in \cite{urli2017constraint} supports the solution of FSM problems over several days.

\begin{algorithm}[H]
\caption{(Subset Algorithm)}\label{algorithm: SA}
\begin{algorithmic}[1]
\State Fix $m\geq 1$
\State Rank the days based on the total load to be delivered on the day, in decreasing order. \label{step: rank}
\State Consider the first $m$ days in the rank and denote them by $i_1,\dots, i_m$.
\State Set the fixed cost to $\frac{b_tm}{|I|}$ \label{step: normalisation}
\State Solve the FSM problem for the first $m$ days and name the so obtained fleet $F$, \label{step: critical step}
\For {$j \in I$ such that $j \not\in \{i_1,\dots,i_m\}}$:
\State adjust $F$ to accommodate vehicles' compatibilities of day $j$ \label{step: comaptibilities}
\State Solve $j$ with $F$
\If {$j$ is infeasible with $F$}
\State Set the routing cost of $j$ to infinity
\EndIf
\EndFor
\State Return $z(m) = bF + \sum_{i \in I}$(routing cost on $i$) \label{step: return}
\end{algorithmic}
\end{algorithm}

The SA solves a ``restricted'' fleet design problem on the first $m$ biggest day, obtaining a fleet $F$. Then it proceeds to solve the routing problem only on all the other days. Note that, for a fixed day $i \notin \{i_1,\dots,i_m\}$ the fleet $F$ is not guaranteed to be feasible. In an attempt to try to prevent infeasibility we adjust $F$ in Step \ref{step: comaptibilities}. This is done only if vehicles' compatibilities imply we need to add some vehicles to F. As an example: if $F$ has no vehicle of type $0$ but there are customers in day $i$ can which only be served by this type of vehicle we add the minimum number of vehicles of type $0$ needed by F. Note that this might not be enough to prevent infeasibility, as for example the placement and time windows of the customers could make a day infeasible even for the adjusted fleet $F$, however there is no trivial way to prevent this possibility. If there is one infeasible day with the given fleet $F$ we simply set the routing cost to infinity. The fixed costs used to solve the restricted fleet design are normalised (Step \ref{step: normalisation}) but then the fixed costs of $F$ are computed using the original vector $b$ (Step \ref{step: return}). Note that there is a trade-off in increasing the number of considered days $m$. For small values of $m$ the produced fleet might be too sub-optimal and probably not feasible on some other day. For higher values of $m$ the fleet has a higher degree of guaranteed coverage but Step \ref{step: critical step} becomes computationally challenging. We run the algorithm SA with $m = 1,\dots,5$ to give SA a fair chance of finding the best possible results it can reach. This way we enable a fair comparison with RMH and BAP. However, except for one run, the best results were found with $m=3$.

For the UF method we run the FSM-solver on each day for $20$ minutes,
which is enough to reach convergence. We used $3$ parallel cores, which amounts to a total of $500$ minutes. To have a fair comparison we run
RMH with a maximum time of $160$ minutes. Once the
execution at the root node is finished we allow an additional $160$
minutes to BAP. We run the SA algorithm, for a fixed $m$, allowing $20m$ minutes to the multi-day FSM-problem solve in Step \ref{step: critical step} and $20$ minutes to every other day. Although we only take the best results of the runs with different values for $m$. This clearly gives SA an advantage (a $5$ times higher computational time), but we will show that, regardless of this advantage, RMH and BAP still outperform SA. Since all methods are stochastic we run each algorithm
$5$ times on each instance and take the average in order to have a
robust result.


\redmark{Given the computational difficulty of the problem, we do not
have a lower bound available with which to compare our
solution. Moreover we do not have the current solution adopted by the
company. However we can compute an ``approximate lower bound'' that
will serve as a lower bound for comparison purposes. It is calculated
assuming we can change the fleet each day. While this bound still
relies on the accuracy of the heuristic technique, it does provide a
conservative bound on the best possible solution that could be
obtained using the given FSM solver. It therefor highlights the
efficacy of the different approaches to choosing a fleet for the
overall problem.}

\redmark{
We run the FSM-solver, with fixed costs set to zero, $5$ times on each
day for $20$ minutes (each run) and take the best run for each day. By
summing the best results of all days we obtain an approximated lower
bound on the operational cost.  We then compute a (formal) lower bound
on the fixed cost by extending the technique explained in the
``\textit{Early termination strategy}'' paragraph of Section
\ref{section: implementation details} to the multi-day problem. This
is straightforward and we do not go into details. The sum of the lower
bounds on operational and fixed costs is taken as the overall
approximated lower bound.  }


\makeatletter
\newcommand{\thickhline}{%
  \noalign {\ifnum 0=`}\fi \hrule height 1pt
  \futurelet \reserved@a \@xhline
}
\makeatother
\newcommand\VRule[1][\arrayrulewidth]{\vrule width #1}
\newcolumntype{?}{!{\vrule width 1pt}}

\begin{table}[h]
\centering
\caption{Computational Results. }
\label{table: performance}
  \begin{tabular}{ l | c | c | c | c | c | c | c |}
      \hline
      \multicolumn{8}{|c|}{instance DC1} \\
      \hline
      \multicolumn{1}{|c|}{method} & cost & $\sigma$ & operational & fixed  & veh & idle & gap    \\
      \hline
      \multicolumn{1}{|c|}{UF} & 785727 & 27030 & 440847 & 344880 & 30.4 & 18.7 & 43.1 \%  \\
      \hline
      \multicolumn{1}{|c|}{SA} & 697739 & 18297 & 413939 & 283800 & 24.6 & 10.7 & 27 \%  \\
      \hline
      \multicolumn{1}{|c|}{RMH} & 605781 & 9866 & 420141 & 185640 & 16.2 & 5.1 & 10.3 \%  \\
      \hline
      \multicolumn{1}{|c|}{BAP} & 578164 & 2835 & 416374 & 161790 & 14.2 & 4.5 & 5.4 \%   \\
      \hline
      \multicolumn{8}{c}{} \\
      \hline
      \multicolumn{8}{|c|}{instance DC2} \\
      \hline
      \multicolumn{1}{|c|}{method} & cost & $\sigma$ & operational & fixed  & veh & idle & gap    \\
      \hline
      \multicolumn{1}{|c|}{UF} &  613004 & 9470 & 364214 & 248790 & 22.2 & 14.0 & 34.2 \%  \\
      \hline
      \multicolumn{1}{|c|}{SA} & 588488 & 9250.4 & 342458 & 246030 & 21.4 & 9.5 & 28.8 \%  \\
      \hline
      \multicolumn{1}{|c|}{RMH} & 499821 & 5727 & 348892 & 150930 & 13.4 & 4.4 & 9.4 \% \\
      \hline
      \multicolumn{1}{|c|}{BAP} & 484543 & 4880 & 340932 & 143611 & 12.8 & 4.1 & 6.0 \% \\
      \hline
  \end{tabular}
\end{table}

In Table \ref{table: performance} we report the average overall cost
(cost), its standard deviation ($\sigma$), a breakdown of the cost in
its components (operational and fixed cost), the average number
of vehicles (veh) in the fleet and the average number of idle vehicles
per day (idle). Lastly we report the average gap with respect to the
approximated lower bound.

It is easy to see that RMH outperforms UF and SA. As expected the main saving
is in the fixed costs. 
As expected, both UF and SA, yield bigger fleets with a higher average
of idle vehicles per day.

The fact that RMH is able to produce a choose a fleet that is only
5.4\% more expensive than a solution where a the best fleet is able to
be chosen each day, indicates the efficacy of the approach described.
As we will show later, the performance of the algorithm depends on the
performance of the FSM-solver. However the two methods, RMH and BAP,
best leverage the performance the underlying FSM-solver, and the small
gaps indicate the solutions are close to the best possible result
obtainable with the given FSM-solver.

It is also interesting to see that the average number of idle vehicles
is even further reduced, in proportion, than the average number of
vehicles. When designing a fleet, companies usually try to have a high
utilisation rate, i.e., the average number of days a vehicle remain
idle should be low. Our approach does not consider this criterion
explicitly in the objective, but it is clear from Table \ref{table:
performance} that the reduction of idle vehicles is a consequence of
the efficiency of the method. Applying BAP leads to an improvement in
the overall cost although the difference is not substantial. It is
possible that the BAP method would add more if we had a bigger problem
due to the fact that forcing the FSM-solver to look at different parts
of the search space might lead to substantial improvements. Overall we
conclude that, for these instances, allowing enough time to the RMH
yields a good solution. Because of this observation, in the next
experiments we will only use RMH.

\subsection{Efficiency of underlying FSM-solver}
One question is how the efficiency of the FSM-solver affects the
overall algorithm. It is clear that better FSM-solvers can provide
better solutions for the single days and lead to a better overall
solution. We argue that the Fleet Size and Mix aspect of the
underlying FSM-solver is quite important. The fixed costs of the
vehicles are the only way of communication between the master problem
and the sub-problems (and hence to the FSM-solver). To demonstrate this
point we solve DC$1$ applying RMH and stopping the FSM-solver at
the $i$-th solution (where $i = 1, \dots, 9$) found in the LNS search
during the solution of a sub-problem. Clearly as $i$ grows the
FSM-solver produces better quality solutions. We set a maximum of a
$120$ sec for each sub-problem.

The results are reported in Table \ref{table: solver-efficiency}. A
consistent improvement is seen in the total and fixed cost if $i$ is
increased from $1$ to $6$. This indicates, as it is shown in
\cite{urli2017constraint}, that the FSM-solver finds good quality
solution rather quickly. After that the cost does not decrease
substantially. The routing costs is lowered as $i$ increases but not
as much as the fixed cost. We can conclude that as the FSM-solver
becomes more efficient we observe a greater savings in fixed costs
rather than in the routing costs.

\begin{table}[]
\centering
\caption{RMH runs with limited solutions for the underlying FSM-solver.}
\label{table: solver-efficiency}
    \begin{tabular}{| l | c | c | c | c | c |}
        \hline
        \multicolumn{1}{|r|}{i =}         & 1 & 2 & 3 & 4 & 5 \\
        \hline
        Fixed cost        & 365550 & 331950 & 331950 & 229200 & 219000 \\
        \hline
        Operational cost  & 428523 & 428969 & 427133 & 426214 & 423149 \\
        \hline
        Total cost        & 794073 & 760919 & 759083 & 655413 & 642149 \\
        \thickhline

        \multicolumn{1}{|r|}{i =}         & 6 & 7 & 8 & 9 & all\\
        \hline
        Fixed cost        & 216150 & 228300 & 183000 & 170850 & 180150 \\
        \hline
        Operational cost  & 425984 & 423592 & 425260 & 423035 & 423716 \\
        \hline
        Total cost        & 642134 & 651892 & 608260 & 593885 & 603866 \\
        \hline
    \end{tabular}

\end{table}


\subsection{Number of days}

\begin{figure}[h]
\centering
\subfloat[\emph{per-day cost}.]
{\includegraphics[width=.45\textwidth]{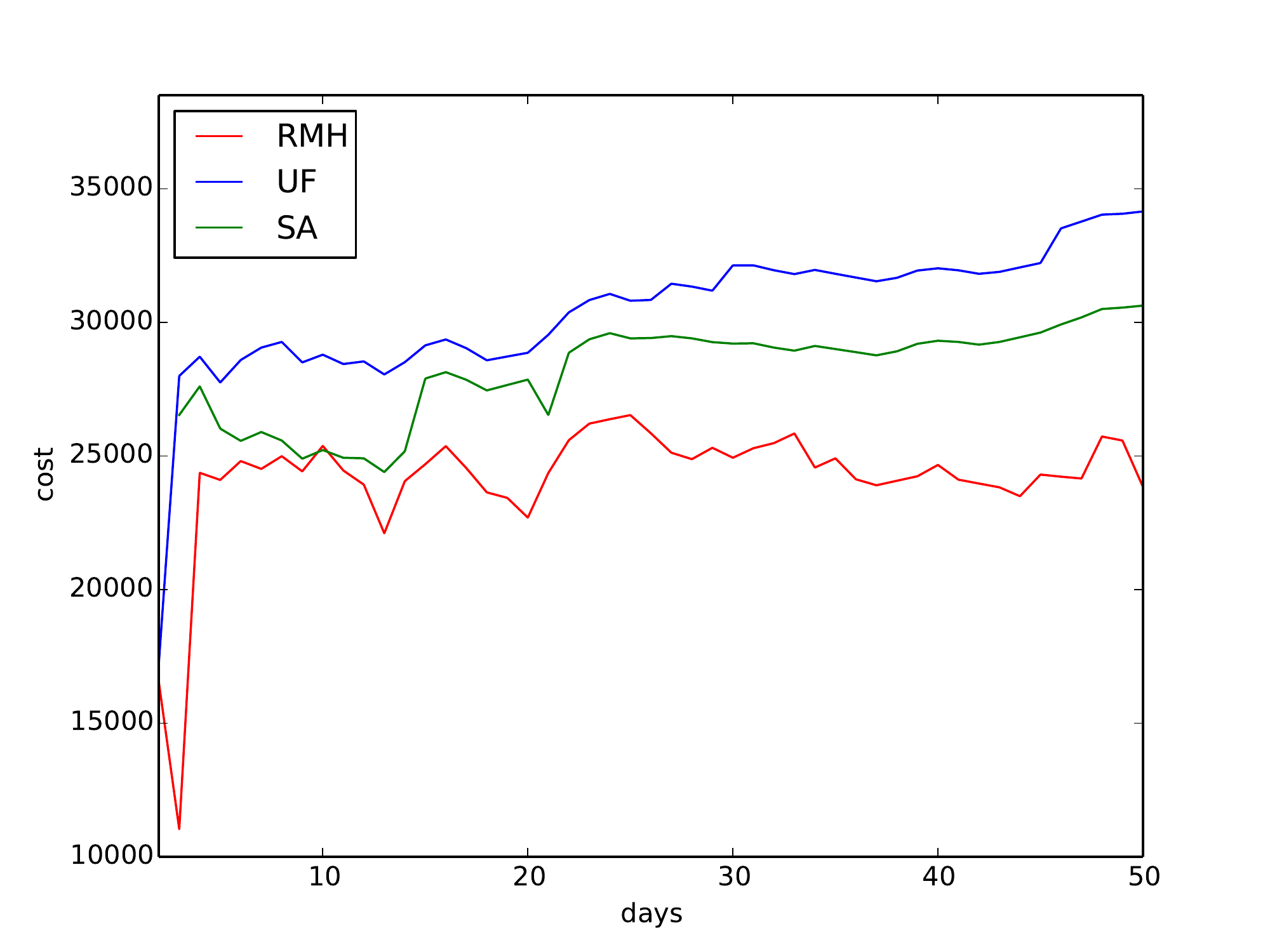}} \quad
\subfloat[\emph{per-day idle vehicles}.]
{\includegraphics[width=.45\textwidth]{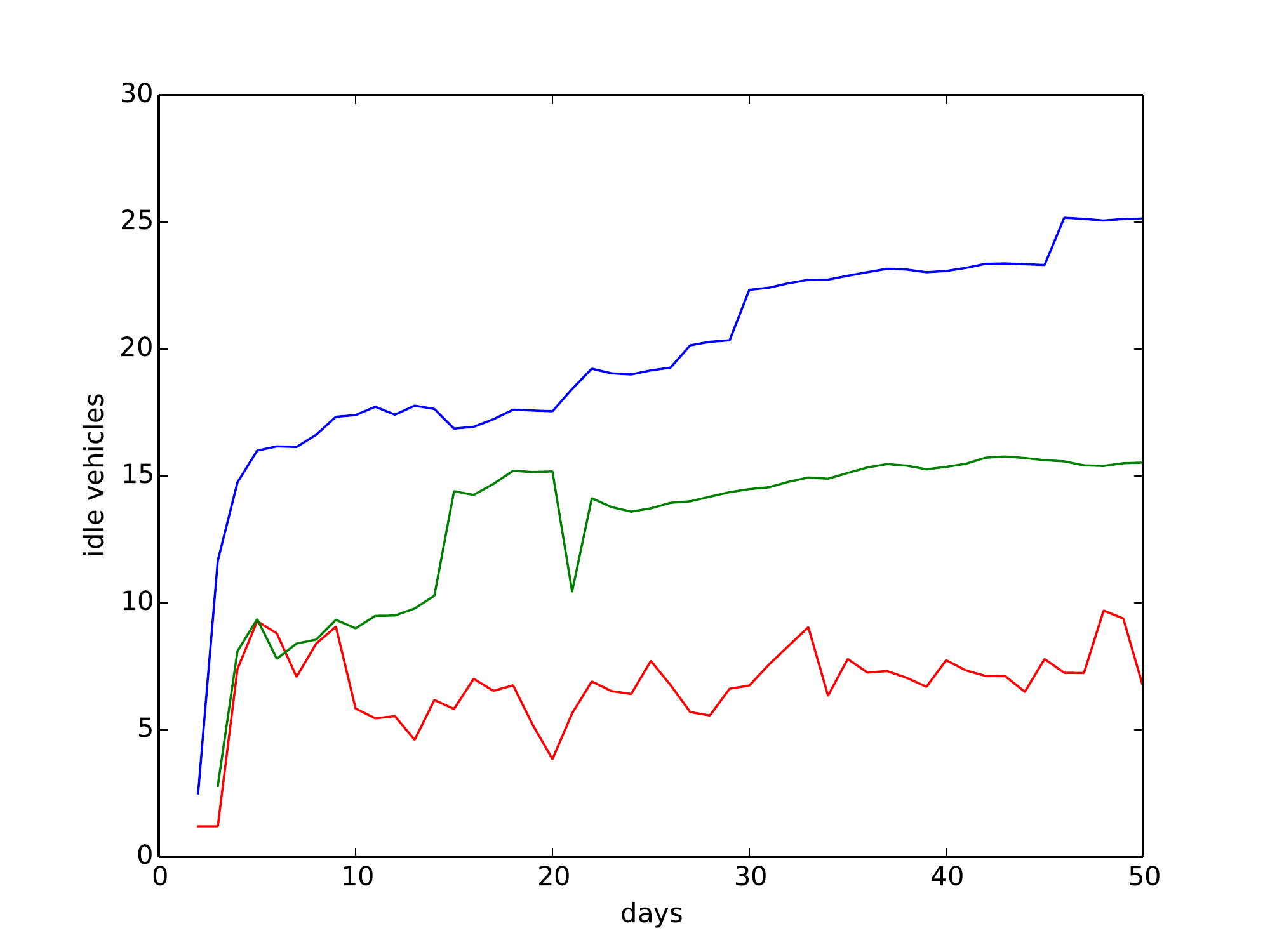}}
\caption{}
\label{figure: day scale}
\end{figure}

The usual meaning of size, when it comes to routing problems, is
the number of customers involved in the daily VRP. This
affects the FSM-solver used to tackle the day sub-problem, rather than
the column generation scheme. The other inputs which can increase the
complexity of the problem are the number of (type of) vehicles and
days. 
We now investigate the performance of RMH, SA and
UF, when the number of days grows. We consider all the $50$ days of demand data
to have a longer horizon available. We denote by $I(d)$ the problem formed by the first $d$ days of demand data. We run UF, SA and RMH on $I(d)$ for $d = 1, \dots, 50$.
We run UF, SA and RMH with time limits as described in section \ref{sec: algorithm performance}. We fixed $m=3$ for SA algorithm as this produced the best results\footnote{In some cases SA produced an infeasible result, i.e., the fleet designed in Step \ref{step: critical step} was not feasible on some other days. In this cases we ran again the whole SA algorithm.}. In Figure \ref{figure: day scale}, for $d = 1, \dots, 50$ we plot (a) the total cost per day and (b) the average number of idle vehicles per day, for problem $I(d)$. It can be seen that RMH is stable when the number of days grow. Oppositely, the fixed cost of the UF tends to deteriorate as the number of days increases. This is due to the fact that, as
the number of days increases, the lack of communication between days
in the solution methods becomes more and more a problem. RMH becomes
more and more effective with respect to UF and SA as we consider more
days. 

The behaviour exhibited by the SA algorithm is harder to interpret. Given all the $50$ days available, the biggest days have indexes $3,21$ and $22$. Therefore once we consider more than $22$ days the fixed cost per day produced by the SA algorithm is constant. The number of days is somehow ignored by the SA, that is the fleet design problem is always solved for $3$ days. Even though this might suggest that SA is not affected by the number of days this is not entirely true. A side effect of adding new days that are not ranked in the first positions is the increase of the average number of idle vehicles, as can be observed in Figure \ref{figure: day scale} (b). Moreover, and more importantly, SA is particularly sensitive to the addition of a new ``hard'' day, i.e. a slight variation in one day might render the obtained fleet infeasible on that day. This can be seen in Figure \ref{figure: day scale} in the area delimited by $d \leq 22$. As mentioned before, after this point the first positions of the rank (Step \ref{step: rank}) do not change. Before this point it is possible to observe how the cost can change rapidly, as the insertions of a new day changes the rank and therefore the fleet. Overall, RMH is more effective, especially when the number of days $d$ grows.

\subsection{Number of vehicle types}

In this section we focus on the second input influencing the complexity of the problem: the number of vehicle types, i.e. $|T|$. In the original DC instance we have $7$ different types of vehicles available. These differ in term of fixed costs, operating costs, capacity, compatibilities with customers and shifts duration. Clearly the last two features play a role in the feasibility of a problem while the first three only influence the cost. We considered instance $DC1$ and removed up to $3$ types of vehicles. Therefore we created a new set of instance having $4$ to $7$ vehicle types available. We call them, respectively $DC1-t, t=4, \dots, 7$. Note that $DC1-7$ and $DC1$ coincide.
We run UF, SA and RMH on each of these instances and report the results in Table \ref{table: vehicle types}. We structure the table so that is possible to observe the behaviour of each method with respect to: the average number of idle vehicles per day, the fixed costs and the operational costs. Clearly a higher number of vehicle types available implies a bigger search space and, coherently, we observe all the three methods worsening as we move from $DC1-5$ to $DC1-7$. However, with less vehicles types, compatibilities constraints have a higher impact on the fleet design. Accordingly, on the instance $DC1-4$ all three methods produced higher costs that on $DC1-4$. The SA method is particularly affected since fleet designed in Step \ref{step: critical step} does not see compatibilities of the other days. Therefore the designed fleet has to be adjusted significantly in Step \ref{step: comaptibilities}, yielding a higher number of idle vehicles.
We can see that both RMH and SA are quite stable. The idle coefficient increases with $|T|$ but not dramatically. On the other side, UF significantly worsen, with respect to the idle vehicles and fixed costs, when we add vehicle types.


\begin{table}[H]
\centering
\caption{Vehicle types.}
\label{table: vehicle types}
    \begin{tabular}{ | l | c | c | c  | }
        \hline
        \multicolumn{4}{|c|}{idle}         \\
        \hline
        instance  &  UF  &  SA  &  RMH        \\
        \hline
        DC1-4     & 8.48 & 9.44 & 5.84        \\
        \hline
        DC1-5     & 8.76 & 12.24 & 3.64       \\
        \hline
        DC1-6     & 16.96 & 12.76 & 4.68      \\
        \hline
        DC1-7     & 18.7  & 10.7 & 5.1        \\
        \hline\hline

        \multicolumn{4}{|c|}{fixed}         \\
        \hline
        instance  &  UF  &  SA  &  RMH      \\
        \hline
        DC1-4     &  208650 &  281550  & 209550   \\
        \hline
        DC1-5     &  204450 &  296400  & 161550   \\
        \hline
        DC1-6     &  301350 & 291150 & 175350     \\
        \hline
        DC1-7     &  344880 & 283800 & 185640     \\
        \hline\hline

        \multicolumn{4}{|c|}{operational}         \\
        \hline
        instance  &  UF  &  SA  &  RMH      \\
        \hline
        DC1-4     &  421623 &  423946 &  425737   \\
        \hline
        DC1-5     &  424626 & 420377  & 420319    \\
        \hline
        DC1-6     &  426314 & 416014  &  418461   \\
        \hline
        DC1-7     &  440847 & 413939  &  420141   \\
        \hline

    \end{tabular}
\end{table}

\section{Concluding Remarks}
\label{section: conclusion}

Strategical and tactical-level decisions in the context of transportation, can
have a significant impact on the  day-to-day operations, and can make the
difference between running a profitable business or not.

In this paper we look at the problem of automatically designing a fleet of
vehicles based on existing historical or forecast demand data. The single-day
version of this problem is known as Fleet Size and Mix (FSM) and has been
extensively studied in literature as an extension of the renowned Vehicle
Routing Problem (VRP). Because the FSM by itself provides little guidance in the
context of a long planning horizon, e.g., one or more years, in this paper we
propose a method to fill this gap. Importantly, our method can leverage any
existing solver for the FSM to solve the more complex problem of designing an
efficient fleet to be operated over a long horizon. The proposed method employs
column generation. Even when the sub-problem is solved heuristically, good
convergence is observed. The method produces either exact or heuristic solutions
to the fleet design problem, depending on whether the sub-problem solver is able
to provide exact or heuristic solutions.

We investigated two frameworks: a restricted master heuristic and a branch \&
price scheme. Both of the methods proposed rely on a suitable FSM-solver. We are
not aware of any other approach able to tackle the same problem in the
literature.

We have compared the technique to the solutions obtained using two common but
simplistic approximations: the union of the best fleets found for each day and the fleet obtained as solution of a critical subset of days.
The proposed method is shown to consistently outperform these heuristics.

While being able to computationally generate an efficient fleet given historical
or forecast demand data is useful, there are several arguments in favour of
exploring techniques to extend an existing fleet and consider the possibility of
hiring, instead of acquiring, vehicles in situations of peak demand. We also
want  to consider options to modify (not necessarily extend) an existing fleet
by selling existing vehicles.

\bibliographystyle{spmpsci}      
\bibliography{paperRefs}   

\begin{thebibliography}{10}
\providecommand{\url}[1]{{#1}}
\providecommand{\urlprefix}{URL }
\expandafter\ifx\csname urlstyle\endcsname\relax
  \providecommand{\doi}[1]{DOI~\discretionary{}{}{}#1}\else
  \providecommand{\doi}{DOI~\discretionary{}{}{}\begingroup
  \urlstyle{rm}\Url}\fi

\bibitem{baldacci2008routing}
Baldacci, R., Battarra, M., Vigo, D.: Routing a heterogeneous fleet of
  vehicles.
\newblock In: The vehicle routing problem: latest advances and new challenges,
  pp. 3--27. Springer (2008)

\bibitem{boschetti2009matheuristics}
Boschetti, M.A., Maniezzo, V., Roffilli, M., R{\"o}hler, A.B.: Matheuristics:
  optimization, simulation and control.
\newblock In: International Workshop on Hybrid Metaheuristics, pp. 171--177.
  Springer (2009)

\bibitem{caceres2015rich}
Caceres-Cruz, J., Arias, P., Guimarans, D., Riera, D., Juan, A.A.: Rich vehicle
  routing problem: Survey.
\newblock ACM Computing Surveys (CSUR) \textbf{47}(2), 32 (2015)

\bibitem{coelho2013thirty}
Coelho, L.C., Cordeau, J.F., Laporte, G.: Thirty years of inventory routing.
\newblock Transportation Science \textbf{48}(1), 1--19 (2013)

\bibitem{crainic2003long}
Crainic, T.G.: Long-haul freight transportation.
\newblock In: Handbook of transportation science, pp. 451--516. Springer US
  (2003)

\bibitem{desrosiers1984routing}
Desrosiers, J., Soumis, F., Desrochers, M.: Routing with time windows by column
  generation.
\newblock Networks \textbf{14}(4), 545--565 (1984)

\bibitem{drexl2012rich}
Drexl, M.: Rich vehicle routing in theory and practice.
\newblock Logistics Research \textbf{5}(1-2), 47--63 (2012)

\bibitem{francis2008period}
Francis, P.M., Smilowitz, K.R., Tzur, M.: The period vehicle routing problem
  and its extensions.
\newblock In: The vehicle routing problem: latest advances and new challenges,
  pp. 73--102. Springer (2008)

\bibitem{gamache1999column}
Gamache, M., Soumis, F., Marquis, G., Desrosiers, J.: A column generation
  approach for large-scale aircrew rostering problems.
\newblock Operations research \textbf{47}(2), 247--263 (1999)

\bibitem{hasle2007industrial}
Hasle, G., Kloster, O.: Industrial vehicle routing.
\newblock In: Geometric modelling, numerical simulation, and optimization, pp.
  397--435. Springer (2007)

\bibitem{hoff2010industrial}
Hoff, A., Andersson, H., Christiansen, M., Hasle, G., L{\o}kketangen, A.:
  Industrial aspects and literature survey: Fleet composition and routing.
\newblock Computers \& Operations Research \textbf{37}(12), 2041--2061 (2010)

\bibitem{joncour2010column}
Joncour, C., Michel, S., Sadykov, R., Sverdlov, D., Vanderbeck, F.: Column
  generation based primal heuristics.
\newblock Electronic Notes in Discrete Mathematics \textbf{36}, 695--702 (2010)

\bibitem{kilby2016fleet}
Kilby, P., Urli, T.: Fleet design optimisation from historical data using
  constraint programming and large neighbourhood search.
\newblock Constraints pp. 1--20 (2016)

\bibitem{kopfer2009combining}
Kopfer, H., Wang, X.: Combining vehicle routing with forwarding: extension of
  the vehicle routing problem by different types of sub-contraction.
\newblock Journal of Korean Institute of Industrial Engineers \textbf{35}(1),
  1--14 (2009)

\bibitem{lubbecke2005selected}
L{\"u}bbecke, M.E., Desrosiers, J.: Selected topics in column generation.
\newblock Operations Research \textbf{53}(6), 1007--1023 (2005)

\bibitem{Gend14Stochastic}
M.~Gendreau O.~Jabali, W.R.: Stochastic Vehicle Routing Problems, chap.~8.
\newblock MOS-SIAM Series on Optimzation. SIAM (2014)

\bibitem{mason1998nested}
Mason, A.J., Smith, M.C.: A nested column generator for solving rostering
  problems with integer programming.
\newblock In: International conference on optimisation: techniques and
  applications, pp. 827--834. Citeseer (1998)

\bibitem{ropke2006adaptive}
Ropke, S., Pisinger, D.: An adaptive large neighborhood search heuristic for
  the pickup and delivery problem with time windows.
\newblock Transportation science \textbf{40}(4), 455--472 (2006)

\bibitem{rossi2006handbook}
Rossi, F., Van~Beek, P., Walsh, T.: Handbook of constraint programming.
\newblock Elsevier (2006)

\bibitem{salhi1993incorporating}
Salhi, S., Rand, G.K.: Incorporating vehicle routing into the vehicle fleet
  composition problem.
\newblock European Journal of Operational Research \textbf{66}(3), 313--330
  (1993)

\bibitem{shaw1998using}
Shaw, P.: Using constraint programming and local search methods to solve
  vehicle routing problems.
\newblock In: International Conference on Principles and Practice of Constraint
  Programming, pp. 417--431. Springer (1998)

\bibitem{urli2017constraint}
Urli, T., Kilby, P.: Constraint-based fleet design optimisation for
  multi-compartment split-delivery rich vehicle routing.
\newblock In: International Conference on Principles and Practice of Constraint
  Programming, pp. 414--430. Springer (2017)

\bibitem{yoshizaki2009scatter}
Yoshizaki, H.T.Y., et~al.: Scatter search for a real-life heterogeneous fleet
  vehicle routing problem with time windows and split deliveries in brazil.
\newblock European Journal of Operational Research \textbf{199}(3), 750--758
  (2009)

\end{thebibliography}

\end{document}